# A Novel Semantics for Belief, Knowledge and Psychological Alethic Modality

**Jonathan J. Mize**

**Abstract** Recently there have been numerous proposed solutions to the problem of logical omniscience in doxastic and epistemic logic. Though these solutions display an impressive breadth of subtlety and motivation, the crux of these approaches seems to have a common theme—minor revisions around the ubiquitous Kripke semantics-rooted approach. In addition, the psychological mechanisms at work in and around both belief and knowledge have been left largely untouched. In this paper, we cut straight to the core of the problem of logical omniscience, taking a psychologically-rooted approach, taking as bedrock the "quanta" of given *percepts*, *qualia* and *cognitions*, terming our approach "PQG logic", short for percept, qualia, cognition logic. Building atop these quanta, we reach a novel semantics of belief, knowledge, in addition to a semantics for *psychological* necessity and possibility. With these notions we are well-equipped to not only address the problem of logical omniscience but to more deeply investigate the psychical-logical nature of belief and knowledge.

---

jonathanmize@my.unt.edu

Department of Philosophy & Religion, University of North Texas, 1155 Union Cir, Denton, TX 76203

## 1 Introduction

Although there has been progress beyond the standard possible-worlds model of belief[1], we have yet to see any substantial progress towards a model of belief that does not lean upon Kripke semantics. Residing in a niche which ostensibly should be fueled by psychologically-geared logical investigations, the fields of doxastic and epistemic logic have instead immersed themselves in continuous refinement of possible-worlds approaches. Nonetheless, we cannot be overly harsh, as there are a few quite convincing, core reasons as to why Kripke semantics have continued to pervade doxastic and epistemic logics. Firstly, after the landmark result of Gödel in 1931, the logic community has all but given up hope of using logic to model such “sticky” issues as belief, desire and emotion. Secondly, though it of course is the case that logic has been used quite extensively in fields such as artificial intelligence, in the modeling of beliefs, desires and emotions and their accompanying structure[2], these approaches have been nearly exclusively of either a probabilistic/“fuzzy” or a possible/impossible worlds nature. As regards the intermingling between psychology and logic, we need only glimpse this quote from the researcher Ben Goertzel’s *Chaotic Logic* (1994) to grok the main idea:

> The early experimental psychologists purposely avoided explaining intelligence in terms of logic. Mental phenomena were analyzed in terms of images, associations, sensations, and so forth. And […] the early logicians moved further and further each decade toward considering logical operations as distinct from psychological operations. It was increasingly realized on both sides that the formulas of propositional logic have little connection with emotional, intuitive, ordinary everyday thought […] The question is whether mathematical logic is a very special kind of mental process, or whether, on the other hand, it is closely connected with everyday thought processes. And, beginning around a century ago, both logicians and psychologists have overwhelmingly voted for the former answer. (p. 42)

Then, later, in the mid-20th century, shortly after the establishment of Kripke semantics in 1959, we witnessed the genesis of the field of epistemic logic. Jaakko Hintikka, in his *Knowledge and Belief*: *An Introduction to the Logic of the two Notions* (1962), presented a manner of interpreting epistemic notions in terms of possible world semantics, laying the groundwork for modern

1 See, most notably and recently, Bjerring & Skipper (2018) and Hawke, Özgün & Berto (2019).

2 For instance, see the BDI (belief-desire-intention) software model and the FLAME (Fuzzy Logic Adaptive Model of Emotions) model. See Kowalczuk & Czubenko (2016) for a comprehensive overview of computational approaches to the modeling of emotion.

epistemic logic. In an insightful piece[3], Hendricks and Rendsvig make the salient point that Hintikka, in his generation of epistemic logic, eschewed the question of what it *is* to know, i.e., the necessary and sufficient conditions of knowledge itself, in favor of pragmatically pursuing the question of what it *means* to know; i.e., the general behavior of knowledge as succinctly and abstracted defined and its logical and informatic permutations and undulations. In [17] all that matters for the semantics of knowledge is "the space of all possible scenarios […] [those] with what I know and those that are incompatible with my knowledge." Now, this is more than adequate if, like stated in [Hendricks, Rendsvig], we are seeking to logically analyze *solely* what it *means* to know. So what then if we desire to also analyze what it *is* to know? For this, we will need to approach belief and knowledge from a more psychological, even metaphysical, standpoint.

Fundamentally, the human experience comprises perceptions, qualia and cognition; "beneath" a belief lies various bundles of these. If we are to fashion a doxastic and epistemic logic geared towards the "what is", as opposed to the "what it means", we must first figure out how perceptions, qualia and cognition intermingle to form knowledge and belief. An ambitious approach could be to examine various logical constructions of and around *emotions*, searching for invariants that may be applicable to doxastic and epistemic logics. Unfortunately, one would find relatively little that is not either firmly buttressed by Kripke semantics [1, 2, 12, 33] or of a fuzzy/statistical [19, 31] nature. Why is this a problem for us? Because—in taking a penetrating, more psychologically-descriptive approach towards belief and knowledge—Kripke semantics and fuzzy logics, these logical frameworks themselves, fail to say anything about how belief and knowledge "hold" in the mind. The following may be a platitude, but it is extremely important to note nonetheless: the mind, that which forms and holds beliefs, relies neither upon possible-worlds semantics nor fuzzy logic; these aforementioned are simply marvelously useful *tools*. What we intend to do here is bring the "tool" closer in line with the true workings of the mind.

## 2 Laying the Foundation

We said above that experience is fundamentally composed of perceptions, qualia and cognition. Granted that we cannot rely upon existing logical frameworks of emotion to sort out the balance of these ingredients, we must seek inspiration elsewhere. Metaphysics, though an oft-discounted source of theoretical insight, is an important source of motivation. In [21] it is pointed out that "quantization" is not merely a process of signal processing or a component of its eponymous subdiscipline of physics; quantization can be *formal* too. From [21], after noting that formal quantization is a well-recognized mathematical necessity: "Sets are quantized as elements, topological spaces are quantized as points, geometry is quantized as lines, which are quantized as points and units of length, and angles, which are quantized as radians or degrees. More generally, any kind of formal system is quantized as symbols representing *objects*, *relationships*, *functions*, and *operations*." (emphasis added)

In first-order logic, objects can be, of course, *variables* or *constants*. And, we have ostensibly all one should need for the relationships ($R_x$), functions ($f_x$) and operations[4] ($\wedge$, $\vee$, $\rightarrow$, etc.). However, one question we must ask is whether an investigation into what belief *is* and how best to logically represent this nature will require us to alter these formal "quanta".

In order to dig deeper into the psychical-logical nature of belief and knowledge, it would behoove us to consider the following formal quanta: *percepts*, *qualia* and *cognitions*. We will consider these, depending on the context, as both variables and constants. For example: ($p_1$, $q_1$, $g_1$) and ($p_x$, $q_x$, $g_x$). More broadly, we will take strings of these quanta to be given "outputs" of higher-level processes. Borrowing some metaphysical inspiration once more from [20, 21], instead of utilizing the traditional state or "instant"-based approach as first seen in Prior's [26], we will adopt what we will call a "meta-simultaneous" approach, having two simultaneous instants, labeled $s^l$ for *linear state* and $s^s$ for *simultaneous state*.

For our $s^s$ states, we will assign to each a set $\mathfrak{F}_{s^s}$ of *volitionary functions*. We may roughly see this as a kind of formal implementation of such theories of rationality and "picoeconomics" as seen in [8], though we will not stress any non-formal construals. Per the set $\mathfrak{F}_{s^s}$, we will have a distinguished function, $f_v$, the *prime volitionary function*, which will take all other *outputs* of $f_n$ of $\mathfrak{F}_{s^s}$ as its *arguments*. In this manner, "intentionality" can be seen as the behavior the function $f_v$, varying in its constituents and its "output" content.

A deceptively similar approach can be seen in [6, 18, 33], where the so-called "STIT logic" is utilized. Many varieties of STIT utilize branching time, shown in [25] to have been originally proposed by Kripke, in a 1958 letter to Prior and subsequently developed by the later, most notably in [26]. In [18], we are given the following definition of a *stit frame*:

$$\langle Tree, <, Agent, Choice \rangle$$

*Choice* is defined as a function mapping each agent $\alpha$ and moment $m$ into a partition $Choice_{\alpha}^{m}$ of the set $H^{m}$ of histories through $m$. Succinctly and naively, each moment has a collection of factors "upon" which the agent can jump from said moment to the next, while the given histories are contingencies of the choice cell $K$ and thus outside $\alpha$'s control.

While STIT logic represents a great technical innovation, when he reflect once more upon the distinction brought up by Hendricks and Rendsvig in [14], we can see that STIT logic merely provides us with more refined methods of pursuing the question of what it *means* to know, as opposed to what it *is* to know. To see this more fully, we need only to consider the semantics of the epistemic operator $\mathsf{K}_{\alpha}$ put forth in [18]:

> Evaluation rule: $\mathsf{K}_{\alpha} A$. *Where $\alpha$ is an agent and m/h an index from an epistemic stit model $\mathscr{M}$,*
>
> - *$\mathscr{M}$, m/h ⊨ $\mathsf{K}_{\alpha} A$ iff $\mathscr{M}$, m′/h′ ⊨ A for all m′/h′ such that m′/h′ $\sim_{\alpha}$ m/h.*
>   - *where m/h $\sim_{\alpha}$ m′/h′ is taken to mean that nothing $\alpha$ knows distinguishes m/h from m′/h′, or that m/h and m′/h′ are epistemically indistinguishable by $\alpha$.*

Although the STIT semantics allow us to deal with indeterminate "timelines" and agent intentionality, it is essentially, for our purposes, Kripke semantics on wheels. That is, this semantics simply transitions from static worlds and accessibility relations to temporal moments and choice cell ($K$) | history ($h$) couplings.

What we have in mind—building from our ($p_1$, $q_1$, $g_1$) and ($p_x$, $q_x$, $g_x$) constants and variables, in addition to our set $\mathfrak{F}_{s}^{s}$ of *volitionary functions*—is a more finely-grained, intra-agentive approach. Belief and knowledge, instead of being defined by external relation (accessibility, choice cell/history coupling), will be defined "internally" as an emergent aspect of the interplay between *percept*, *qualia, cognition* and *volitionary function*. But, before we dive into the technical details, we need to gather a motivating intuition for what exactly belief is.

First and foremost, what we desire to do, is fuse *functionalism* and *dispositionalism*; we want to be sufficiently functionalist so as to provide a rough working outline of psychological mechanism, and we want to be sufficiently dispositionalist so as to form a logical model that can consistently "recognize" those quanta that are relevant to belief. In terms of the criteria of sentential complements of the form *that P*, say, from the sentence *α believes that P,* [31] notes the point that the objects of belief need be individuated very finely, stating, "for almost any two distinct sentential clauses, *that P* and *that Q*, one can find a context where it seems plausible to say that someone believes that *P*, but disbelieves that *Q*." Utilizing quanta of the form $(p_1, q_1, g_1)$ and $(p_x, q_x, g_x)$, we can reach a very finely-grained individuation indeed. With this approach, we will be actively "thinking outside of the belief box", a la Schwitzgebel in his [30].

# 3 The Approach

**Definition 1** (**PQG Model**) Let $W$ be a non-empty set of possible worlds, let $\mathbb{S}_l$ be a set of "linear" moments and let $\mathbb{S}_s$ be a set of "simultaneous" moments[5], $\mathbb{B}$ is a set of *belief states*; $\mathbb{D}$ is a set of *volitionary determination*, per a given $b^{sx}{}_x$, $\mathbb{R}$ is a set of *rules*, per a given $f_v{}^{s_x}$, $\mathfrak{F}_{s^s}$ is a set of *volitionary functions*, $\mathbb{C}$ is a set of *psychological concepts*, $R_{\mathbb{S}}$ is a relation $\langle \mathbb{S}_l, \mathbb{S}_s \rangle \times W$; $R_{\mathbb{B}}$ is a relation $\mathbb{S}_s \times \mathbb{B}$; $R_{\mathbb{C}}$ is a relation $\langle \mathbb{S}_l, \mathbb{S}_s \rangle \times \mathbb{C}$; $R_{\mathbb{D}}$ is a relation $\langle \mathbb{S}_s, b^{sx}{}_x \rangle \times \mathbb{D}$; $R_{\mathbb{R}}$ is a relation $\mathbb{S}_s \times \mathbb{R}$; $R_{\mathfrak{F}}$ is a relation $\mathbb{S}_s \times \mathfrak{F}_{s^s}$; $\supset_s$ is the relation of "containment" of a given $s^l$ within a given $s^s$. A PQG model for a single agent is a structure:

$$\langle W, R, \mathbb{S}_l, \mathbb{S}_s, \mathbb{B}, \mathbb{C}, \mathbb{D}, \mathbb{R}, \mathfrak{F}_{s^s}, R_{\mathbb{S}}, R_{\mathbb{B}}, R_{\mathbb{C}}, R_{\mathbb{D}}, R_{\mathbb{R}}\ R_{\mathfrak{F}}, <_l, <_s, \supset_s, <_h{}^{bs_x}, V \rangle,$$

where we take the ordering relations to be those of a traditional, instant-based model of time and $<_h{}^{bs_x}$ is a so-called called "hypothetical time" relation, to be defined later.[6]

After downplaying the efficiency of possible-worlds semantics above, perhaps we should explain why we chose to incorporate the semantics. Granted that we are developing a semantics of psychologically-geared modality—and, shortly, will present a notion of psychological *alethic* modality—we can still utilize possible-world semantics as a sort of "covering" for our semantics. This holds to the intuition that psychological modality is dependent upon *metaphysical* (possible-worlds) modality. Let us now explain some of our new terminology and how it relates to our necessarily higher-order language $\mathcal{L}$.

We said previously that we wish to take both a functionalist and dispositionalist approach to belief. With the belief state, $b_n \in \mathbb{B}$, we utilize the dispositionalist portion. Instead of taking a given belief as representation with content, the approach criticized in [30], we construe a given

5 Although we have neglected to discuss sorts for concision, we will take $\mathbb{S}_l$ and $\mathbb{S}_s$ as *sorted* sets of moments. For an example of a many-sorted modal logic, see Leuştean, Moangă & Şerbănuţă (2018).

6 Our temporal relations, $<_l$, $<_s$ and $<_h{}^{bsx}$ will also be *sorted*. For an example of a many-sorted temporal logic, see Enjalbert & Michel (1984).

belief state, $b_n$, as a hypothetical set of quanta and functions. In essence, a belief state is the resulting "bundle" of quanta, were the object of a given belief to be actualized, i.e., availed to the agent in question. Belief states are paired with given moments, $s^s$, by the relation $R_{\mathbb{B}}$. Thus we have the (preliminary) definition of a random belief state:

**Definition 2** (**Belief State**) $$b^{sx}{}_3 = \{p_2 \mapsto g_3 \mapsto q_1\},$$

where the functions $\mapsto$ are regarded endogenously per the given $s^s$; where we consider a set of quanta, $Q$, we can either consider it with or without the functions $\mapsto$, denoted by, respectively, simply $Q^{\mapsto}$ and $Q$.

Next, before we define members of $\mathbb{C}$, what we call *psychological concepts*, we will define two functions, $\mathcal{T}$, a *taking function*, and $\mathcal{F}$ a *forming function*.

**Definition 3** (**Taking Function**) Let $Q^{\mapsto}{}_x$ and $Q^{\mapsto}{}_y$ be given sets of quanta, such that $Q^{\mapsto}{}_x >_l Q^{\mapsto}{}_y$, then let $\mathcal{T}$ be defined as $\mathcal{T}$: $Q^{\mapsto}{}_x \mapsto Q^{\mapsto}{}_y$. Note that it is possible that the domain of $\mathcal{T}$ is $Q^{\mapsto} \subset Q^{\mapsto}{}_x$ or $Q \subset Q^{\mapsto}{}_x$.

**Definition 4** (**Forming Function**) Let $Q^{\mapsto}{}_y$ be the result of an application of $\mathcal{T}$ on a given $Q^{\mapsto}{}_x$, and let $Q^{\mapsto}{}_z$ be any given set of quanta, then let $\mathcal{F}$ be defined as $\mathcal{F}$: $Q^{\mapsto}{}_y \mapsto Q^{\mapsto}{}_z$.

We can now define our notion of psychological concepts, the intuition behind which is one of "conceptual comparison", whereby a currently experienced set of quanta is taken back to (through memory) a certain salient set of quanta. From this, the two sets of quanta are compared according to certain salient (higher-order) characteristics and a certain "output" is recommended. A psychological concept is then the precise nature of such comparison to "output" mapping. Psychological concepts can be seen as the building blocks of choice, and, later, belief and knowledge.

**Definition 5** (**Psychological Concepts**) Let $\mathcal{F}$: $Q^{\mapsto}{}_y \mapsto Q^{\mapsto}{}_z$ be any given application of $\mathcal{F}$, we define a concept of the set $\mathbb{C}$, $c_x$, to be a given mapping of the function $\mathcal{F}$, $Q^{\mapsto}{}_y \mapsto Q^{\mapsto}{}_z$. Note that we can have multiple concepts per given $Q^{\mapsto}$ of $s^l$.

Although the following point is rather obvious, we will briefly stress it nonetheless—the language $\mathcal{L}$ of our PQG system must be of a higher-order nature. This can be seen by examining the notion of psychological concepts. These will be fed into the agent's set $\mathfrak{F}_{s^s}$ of volitionary functions, per $s^s$, and will subsequently be quantified over. Now, having explicated $\mathcal{T}$, $\mathcal{F}$ and members of $\mathbb{C}$, we can now define volitionary functions.

**Definition 6 (Volitionary Functions)** Let the contents of $\mathfrak{F}_{s^s}$ be an $n$-tuple, $f_v, f_i, \ldots, f_n$ such that $f_v$, the *prime volitionary function*, takes the *outputs* of $f_i, \ldots, f_n$ as its arguments; let $f_i$ be such that it takes the *arguments* of $f_j, \ldots, f_n$ as its arguments and so on. Per the nature of the arguments they take, inversely, each volitionary function is assigned an order, such that, in the above case, we would have $f_v^{0}$, $f_i^{1}$, $f_j^{2}$, $f_n^{n+1}$. In terms of the arguments of the members of $\mathfrak{F}_{s^s}$, each function will be of the following general form,

$$f_x^n: (\langle c_i, \mathsf{Q}^{\mapsto}{}_z{}^i\rangle, \langle c_j, \mathsf{Q}^{\mapsto}{}_z{}^j\rangle, \ldots, \langle c_n, \mathsf{Q}^{\mapsto}{}_z{}^n\rangle) \mapsto \mathsf{Q}^{\mapsto},$$

where each $\mathsf{Q}^{\mapsto}{}_z$ is from a forming function, $\mathcal{F}$: $\mathsf{Q}^{\mapsto}{}_y \mapsto \mathsf{Q}^{\mapsto}{}_z$ and where the output of $f_x^n$, $\mathsf{Q}^{\mapsto}$ is to be called a *recommended quanta string*, or RQS.

Notice how above we defined $f_v$ as a function taking the *outputs* of $f_i, \ldots, f_n$ as its arguments. Thus, we have,

$$f_v: ((\mathrm{RQS}^{\,f_i}), (\mathrm{RQS}^{\,f_j}), \ldots, (\mathrm{RQS}^{\,f_n})) \mapsto \mathsf{Q}^{\mapsto},$$

where the output $\mathsf{Q}^{\mapsto}$ is any given member $\mathsf{Q}^{\mapsto}{}_{sl}$.

If we were so-inclined to "dress up" our language, we would have a unique preference relation $\prec$ assigned to each volitionary function. However, in this paper, since we will only be focusing on the formation of doxastic and epistemic notions, we will gloss over such details. Having explicated the above, we are now only a few definitions away from implementing our semantics for belief.

**Definition 7 (Pre-Belief State Moments)** Given each $b^{sx}{}_x \in \mathbb{B}$, we *can* (restrictions soon to be defined) have a set $\mathbb{P}$ of *pre-belief state* moments. Where a given $b^{sx}{}_x$ is the output of a given $f_v$: $((\mathrm{RQS}^{\,f_i}), (\mathrm{RQS}^{\,f_j}), \ldots, (\mathrm{RQS}^{\,f_n})) \mapsto \mathsf{Q}^{\mapsto}$, a pre-belief state moment $p^{sx}{}_x$ is a given moment $s^l$ such that $s^l \supset_S s^s$, where the given $f_v$ outputting the belief state moment $b^{sx}{}_x$ (which is just a $\mathsf{Q}^{\mapsto}$) is such that $f_v \supset_S s^s$. Following the fact that beliefs need not be satisfied in the world to be beliefs, i.e., beliefs need not materialize into events of recognition and, thus, *knowledge*, we order members of $\mathbb{P}$ along an "orthogonal" time axis, ordered by the relation $<_h^{bs_x}$, what we will call the *hypothetical time* relation per a given belief state. Members of $\mathbb{P}$ are *not* ordered by $<_l$, although the $s^s$ are still ordered, of course, by $<_s$. This, along with the interaction of $<_s$ with $<_l$ entails that the relation $<_s$ is a "dependent relation", although we will not bother to explicate this in technical detail.

We said previously that—in taking the *dispositionalist* approach to belief—belief state moments represent strings of quanta, *were* the "item" of belief to be introduced. That is, the attainment of a belief state is the *hypothetical* satisfaction of currently held *invariants*. What exactly do we mean by invariants? Simply, the answer is in the name, *disposition*alism. In our PQG logic, we take a disposition to believe to be certain invariant aspects of an agent's prime volitionary function, $f_v$. In order to best capture the notion of volitionary invariance, we develop a distinguished set $\mathbb{D}_{fv}$, the *set of volitionary determination* for a $f_v$ at a given $s^s$ moment.

**Definition 8 (Volitionary Determination)** Take any given $f_v^{s_x}$, per an $R_{\mathbb{B}}$ determined $b^{s_x}{}_x$ at $f_v^{s_x}$'s moment $s^s$, we define the set $\mathbb{D}_{fv}{}^{b_x}$, the *set of volitionary determination* for agent *α*, per the given belief state, to be a set of "rules" of the nature of all $f_x{}^n$ such that $f_v^{s_x}$ takes all outputs of all $f_x{}^n$ as its arguments. These are seen to be rules dictating the "ability" of $f_v^{s_x}$ to "reach" $b^{s_x}{}_x$, such that, if these rules hold, per all given $p^{s_x}{}_x$, then $b^{s_x}{}_x$ will eventually hold. We also permit a *subset* of the set $\mathbb{D}_{fv}{}^{b_x}$, a set $\mathfrak{S}_{fv}{}^{b_x}$ of *minimal rules* and a *superset*, a set $\mathfrak{X}_{fv}{}^{b_x}$ of *maximal rules*. We will see shortly that this partitioning helps move us towards our semantics of psychological alethic modality. Each $f_v^{s_x}$ has an accompanying set $\mathbb{R}_{fv}$ of rules along each composite function, such that we may "compare" $\mathbb{R}_{fv}$ and $\mathbb{D}_{fv}{}^{b_x}$. Below we list a hypothetical set $\mathbb{D}_{fv}{}^{b_2}$ consisting of only two functions "below" the prime volitionary function $f_v$,

$$f_i{}^1: (P^i{}_2\langle(c_i), (Q^{\mapsto}{}_z{}^i)\rangle, P^i{}_4\langle(c_j), (Q^{\mapsto}{}_z{}^j)\rangle) \mapsto P^j{}_4(Q^{\mapsto})$$

$$f_j{}^2: (P^i{}_6\langle(c_i), (Q^{\mapsto}{}_z{}^i)\rangle, P^i{}_2\langle(c_j), (Q^{\mapsto}{}_z{}^j)\rangle) \mapsto P^j{}_3(Q^{\mapsto}),$$

where the $P^x{}_x$s are given predicates, whose definitions may be of various sorts. We provide example definitions below,

$$P^i{}_6\langle(c_i), (Q^{\mapsto}{}_z{}^i)\rangle \equiv c_i \ni \forall\, Q^{\mapsto}{}_y\, (Q^{\mapsto}{}_y \in \mathcal{F})\, [Q^{\mapsto}{}_y \in RQS_{s_x}{}^{f_i}] \ldots Q^{\mapsto}{}_z{}^i \ni \exists\, g, q\, ((g, q) \in \mathcal{T}\!:\; Q^{\mapsto}{}_y >_l Q^{\mapsto}{}_x{}^2$$

$$P^j{}_3(Q^{\mapsto}) \equiv Q^{\mapsto} \ni \exists\, s^l\, (s^l <_l Q^{\mapsto})$$

Clearly there is a vast potential for various predicate varieties, and we will save a precise investigation thereof for another date. Nonetheless, from this brief exposition, it should be evident how the definition of $\mathbb{D}_{f_v}{}^{b_x}$ ties into the PQG framework at large. Note that we leave open the possibility of having sets of volitionary determination for generalized (non-doxastic) $\mathsf{Q}^{\mapsto}$ sets, symbolized

**Definition 9** (**Volitionary Acceptance | Volitionary Invariance**) If we have any given $f_v{}^{s_x}$, per an $R_{\mathbb{B}}$ determined $b^{sx}{}_x$, and, if the current set $\mathbb{R}_{f_v}$ is such that $\forall f\,(f \in \mathbb{D}_{f_v}{}^{b_2})\,[f \in \mathbb{R}_{f_v}]$, then $b^{sx}$ is said to have attained *volitionary acceptance* by $f_v{}^{s_x}$ at $s^l, s^{\mathsf{s}} \ni s^l \supset_{\mathsf{s}} s^{\mathsf{s}}$, symbolized like, $b^{sx} \diamond f_v{}^{s_x}$. For each $s^l, s^{\mathsf{s}} \ni s^l \supset_{\mathsf{s}} s^{\mathsf{s}}$, before $b^{sx}{}_x$ is output as its corresponding $\mathsf{Q}^{\mapsto}$, if each $f_v{}^{s_x}$ admits volitionary acceptance, then the sequence of moments, $(s^l, s^{\mathsf{s}})_1, \ldots, (s^l, s^{\mathsf{s}})_n$ leading up to the output $\mathsf{Q}^{\mapsto}$ of $b^{sx}{}_x$, is said to exhibit *volitionary invariance* under the sequence $f_v{}^{s_x}{}_1, \ldots, f_v{}^{s_x}{}_n$, symbolized like, $[(s^l, s^{\mathsf{s}})_1, \ldots, (s^l, s^{\mathsf{s}})_n] \between [f_v{}^{s_x}{}_1, \ldots, f_v{}^{s_x}{}_n]$.

Above we gave a preliminary (though essentially comprehensive) definition of *pre-belief state* moments, members of $\mathbb{P}$, stating that there were certain restrictions. The restriction is that, in order for there to be a pre-belief state or sequence thereof, the given $b^{sx}$ must have attained *volitionary acceptance* by $f_v{}^{s_x}$, or, for a sequence, the moments $(s^l, s^{\mathsf{s}})_1, \ldots, (s^l, s^{\mathsf{s}})_n$ must exhibit *volitionary invariance* under the sequence $f_v{}^{s_x}{}_1, \ldots, f_v{}^{s_x}{}_n$. Otherwise, we have $M, s^{lx} \nvDash b^{sx}{}_x$. For pre-belief state moments, $p^{sx}{}_x$, there is a sequence $p^{sx}{}_1, \ldots, p^{sx}{}_n$, as long as the sequence $f_v{}^{s_x}{}_1, \ldots, f_v{}^{s_x}{}_n$ at all $p^{s_{\mathsf{s}}}{}_x \supset_{\mathsf{s}} p^{s_l}{}_x$ affords volitionary invariance. Pre-belief state moments are simply volitionally invariant $s^l, s^{\mathsf{s}}$.

Before we introduce the doxastic, epistemic and alethic semantics for PQG logic, we need to discuss one last topic, that of nested doxastic and epistemic operators, viz. $\mathfrak{B}\mathfrak{B}p$, $\mathfrak{K}\mathfrak{K}p$. In the standard, possible-worlds semantics of doxastic logic, nested operators are simply evaluated along the extent of the corresponding $R$-accessibility; $\mathfrak{B}\mathfrak{B}p$ and $\mathfrak{K}\mathfrak{K}p$ are $\Box\Box p$ which is simply the evaluation of all $w''$ such that $w''R\,w'$. In PQG logic, as our semantics for belief does not utilize possible-worlds semantics (via its essential *psychologico-descriptive* nature), we necessarily devise a different criteria of satisfaction for $\mathfrak{B}\mathfrak{B}p/\mathfrak{K}\mathfrak{K}p$. We introduce a meta-belief/meta-knowledge operator, $\mathfrak{B}^{\mathfrak{M}}/\mathfrak{K}^{\mathfrak{M}}$ that can ascend an arbitrary number of "degrees", as $\mathfrak{B}^{\mathfrak{M}^n}p/\mathfrak{K}^{\mathfrak{M}^n}p$. Now, we may introduce our definitions of satisfaction.

**Definition 10** (**Satisfaction**) Where our $\mathcal{L}$ is of an arbitrarily higher-order (and, *many-sorted*, but we will omit the specifics for concision)[10], such that we have sets of predicate variables of ascending order, $\{{}^0X_1{}^n, \ldots, {}^0X_n{}^n\}, \ldots, \{{}^nX_1{}^n, \ldots, {}^nX_n{}^n\}$ and function variables of ascending order, $\{{}^0F_1{}^n, \ldots, {}^0F_n{}^n\}, \ldots, \{{}^nF_1{}^n, \ldots, {}^nF_n{}^n\}$. $\mathcal{L}$ is then defined the usual inductive way from a set $\Phi$ of atomic sentences, the quantifiers $\forall$, $\exists$, a set $\{\neg, \wedge, \vee, \rightarrow, \equiv\}$ of connectives, a set $\{x_1, \ldots, x_n\}$ of *individual* variables and a set $\{c_1, \ldots, c_n\}$ of constants.

Given that each $s^l, s^s \ni s^l \supset_s s^s$ is a member of a given $w \in W$, the satisfaction of $\Box p$ and $\Diamond p$ and the usual temporal operators are the same as usual. Although we give definitions of satisfaction in terms of moments $s^{sx}$ and $s^{lx}$, definitions in terms of worlds are equivalent, in that $(M, s^{sx}, s^{lx} \vDash \mathfrak{B}p) \equiv (M, w \vDash \mathfrak{B}p)$ when $(s^{sx}, s^{lx}) \in w$.

For any $(s^{sx}, s^{lx}) \ni s^{lx} \supset_s s^{sx}$:

(B1) $M, s^{sx}, s^{lx} \vDash \mathfrak{B}p$ *iff* $b^{sx} \in s^s$ attains *volitionary acceptance* by $f_v{}^{sx}$ and the sequence of moments, $(s^l, s^s)_1, \ldots, (s^l, s^s)_n$ leading up to the output $\mathrm{Q}^{\mapsto}$ of $s^{lx}$ exhibits *volitionary invariance* under the sequence $f_v{}^{sx}{}_1, \ldots, f_v{}^{sx}{}_n$.

(B1.1) $M, s^{sx}, s^{lx} \vDash \mathfrak{B}(p \rightarrow q)$ *iff* whenever $p$ is true at a pre-belief state moment, so is $q$.

(B1.2) $M, s^{sx}, s^{lx} \vDash \mathfrak{B}(q \rightarrow p)$ *iff* whenever $q$ is true at a pre-belief state moment, so is $p$.

(B1.3) $M, s^{sx}, s^{lx} \vDash \mathfrak{B}(q \equiv p)$ *iff* whenever $q$ is true at a pre-belief state moment, so is $p$, as well as the converse.

(B2) $M, s^{sx}, s^{lx} \vDash (\mathfrak{B}p \rightarrow \mathfrak{B}q)$ *iff* $p$ and $q$ are given belief states of the given simultaneous moment, viz., $(b^{sx}{}_x, b^{sx}{}_y) \in s^{sx}$ and if whenever $p$ holds true that $\vDash \mathrm{Q}^{\mapsto}(b^{sx}{}_x)$, then $q$ holds true as $\vDash \mathrm{Q}^{\mapsto}(b^{sx}{}_y)$.

(B3) $M, s^{sx}, s^{lx} \vDash \mathfrak{B}^{\mathfrak{M}n}p$ *iff* conditions for $\vDash \mathfrak{B}p$ are met, and then, <u>while</u> a given $f_{vi}, \ldots, f_{vn}$ $(\mathbb{D}_{f_v}{}^{s_{lx}})$ holds such that $b_{\mathfrak{M}n}{}^{sx}$ attains volitionary acceptance by $f_{vi}, \ldots, f_{vn}$; $\mathfrak{B}^{\mathfrak{M}n}p$ is then said to hold as long as the moments $(s^l, s^s)_1, \ldots, (s^l, s^s)_n$ leading up to the output $\mathrm{Q}^{\mapsto}$ of $b_{\mathfrak{M}n}{}^{sx}$ exhibit *volitionary invariance* under the sequence $f_{vi}, \ldots, f_{vn}$ of $f_{vi}, \ldots, f_{vn}$ $(\mathbb{D}_{f_v}{}^{s_{lx}})$.

(B4) $M, s^{s_x}, s^{l_x} \vDash \mathfrak{K}p$ *iff* conditions for $\vDash \mathfrak{B}p$ are met, <u>while</u> for the given $s^{l_x}$ or, equivalently, for all $s^{l_n}$ in question, $M, s^{l_n} \vDash p$; that is, $p$ is a given $\vDash \mathsf{Q}^{\mapsto}$ such that $p \notin p^{sx}{}_x$.

(B5) $M, s^{s_x}, s^{l_x} \vDash \mathfrak{K}^{\mathfrak{M}n}p$ *iff* conditions for $\vDash \mathfrak{B}p$ <u>and</u> $\vDash \mathfrak{B}^{\mathfrak{M}n}p$ are met, <u>while</u> for the given $s_l^x$ or, equivalently, for all $s_l^n$ in question, $M, s_l^n \vDash p$ <u>and</u>, for the moments $(s^l, s^s)_1, \ldots, (s^l, s^s)_n$ leading up to the output $\mathsf{Q}^{\mapsto}$ of $b_{\mathfrak{M}n}{}^{sx}$ $M, s_l^n \vDash p$.

For any $s_l^x$:

(S1) $M, s^{l_x} \vDash \boxdot p$ *iff* <u>all</u> members of the superset $\mathfrak{X}_{fv}{}^{b_x}$ of maximal rules of the corresponding $\mathbb{D}_{fv}{}^{s_{lx}}$ are satisfied and the sequence of moments, $(s^l, s^s)_1, \ldots, (s^l, s^s)_n$ leading up to the output $\mathsf{Q}^{\mapsto}$ of $s^{l_x}$ exhibits *volitionary invariance* under the sequence $f_v{}^{s_x}{}_1, \ldots, f_v{}^{s_x}{}_n$.

(S2) $M, s^{l_x} \vDash \circledcirc p$ *iff* <u>not all</u> members of the of the corresponding $\mathbb{D}_{fv}{}^{b_x}$ are satisfied (while all members of the subset $\mathfrak{S}_{fv}{}^{b_x}$ of minimal rules are), but $b^{sx}$ attains *volitionary acceptance* by $f_v{}^{s_x}$ and the sequence of moments, $(s^l, s^s)_1, \ldots, (s^l, s^s)_n$ leading up to the output $\mathsf{Q}^{\mapsto}$ of $s^{l_x}$ <u>does not</u> exhibit *volitionary invariance* under the sequence $f_v{}^{s_x}{}_1, \ldots, f_v{}^{s_x}{}_n$.

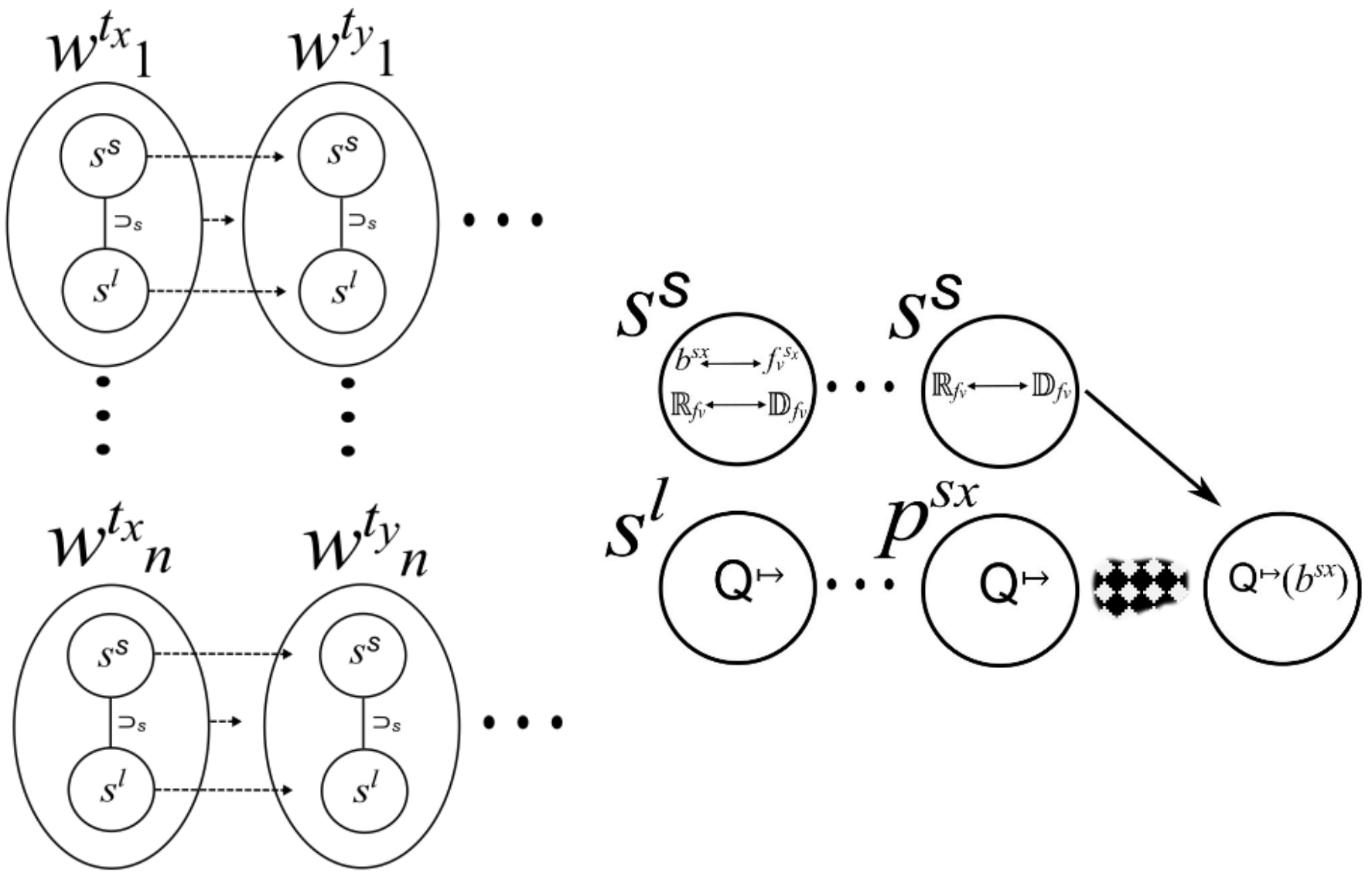

**Fig. 1** On the left we see an illustration of given worlds $\{w^{t_x}{}_1, \ldots, w^{t_x}{}_n\}$ and $\{w^{t_y}{}_1, \ldots, w^{t_y}{}_n\}$, per given times $t_x$ and $t_y$. We can see the "containment" of $s^s$ and $s^l$ within each world, in addition to the "containment" of $s^l$ within $s^s$. On the right we see an illustration of the "contents" of paired moments, *linear*, *simultaneous* and *pre-belief state*. At $s^s$

we see a given belief state attaining volitionary acceptance by a given volitionary function, and we also see a given set of rules “measuring up” to a set of volitionary determination per the belief state. The figure on the right contains all the necessary ingredients for $M, b^{sx} \in s^{l_x} \models \mathfrak{B}p$. The polka-dotted “patch” we see between the right most $\mathsf{Q}^{\mapsto}$ and the $\mathsf{Q}^{\mapsto}( b^{sx}{}_{x})$ represents the “jump” from hypothetical ordering, $<_h^{bsx}$, to regular, linear ordering $<_l$.

# 4 Consideration of Axioms

It is clear that—given that we opted to build a system geared towards the “what is” of belief and knowledge, as opposed to the “what it means”—our axioms will be of a peculiar nature. However, it will be an interesting endeavor to see how our potential *psychologico-descriptive* axioms compare to various varieties of the traditional *exogenously-reasoned* axioms. Before we commence our investigation however, let us visit a source upon which we can further cement our previous distinction between “psychologico-descriptive” axioms and “exogenously-reasoned” axioms. Here, in [27], we see a quite perceptive—especially granted the date of publication—analysis of the relation of axiomatic systems to psychology.

> [C]an these fundamental axioms [of formal logic] be considered practical precepts based on psychological laws ? If so, *what are these fundamental psychological laws* ? If they are not distinguishable from the logical axioms, and these last are therefore laws of nature, *how are the fallacies which consist in their violation possible*? The distinction between nature or 'things' and our thinking about things, will hardly help us here, for these axioms of logic are at once statements about things and about the necessities of our thought. Here, then, we are face to face with a difficulty which is just one aspect of the problem, ' How is knowledge possible ?' with its companion problem, 'How is error possible ?' (emphasis added)

Here, in the italicized portion most specifically, we see the crux of the issue we have been getting at. We establish a distinction between (i) the functional assembly of cognition; the entire logically analyzable apparatus from perception and memory to reasoning and qualia and (ii) an apparatus of implication, with all of the accompanying functional assembly abstracted away. The former is what PQG logic endeavors to utilize and analyze, and the latter is what possible-worlds semantics affords.

Although it may initially seem to be a useless “sacrifice of expressivity”, the genesis of a logic upon the former half of the distinction represents a necessary and exciting foray into uncharted territory. Now, we shall examine the axiomatics.

In [23] we see a survey of the field of proposed epistemic logics. Here is the overview presented, “Proposed logics have included **KT** (Williamson2000, Sect.10.4), **S4** (Hintikka1962), **S4.2** (Lenzen1978; Stalnaker2006), **S4.3** (van der Hoek1996), **S4.4** (von Kutschera1976) and **S5** (Fagin et al.1995).” Instead of erecting a PQG-tailored set of conditions (due to space constraints), we will simply compare and contrast PQG with each logic above. Below we will commence a novel notation, using $\mathfrak{D}^{\mathfrak{B}}\varphi$ instead of $\mathfrak{B}\varphi$ and $\mathfrak{D}^{\mathfrak{K}}\varphi$ instead of $\mathfrak{K}\varphi$, highlighting the difference between PQG and traditional doxastic and epistemic logics; where $\mathfrak{D}$ stands for *descriptive*.

### 4.1 PQG vs. KT

Starting with the frame condition for **K**,

$$K(\varphi \to \psi) \to (K\varphi \to K\psi),$$

$$\mathfrak{D}^{\mathfrak{K}}(\varphi \to \psi) \to (\mathfrak{D}^{\mathfrak{K}}\varphi \to \mathfrak{D}^{\mathfrak{K}}\psi),$$

given the satisfaction definitions for PQG, we have,

- $M, s^{sx}, s^{l_x} \vDash \mathfrak{D}^{\mathfrak{K}}(\varphi \to \psi)$ *iff* $\varphi$ and $\psi$ are at given pre-belief state moments and whenever $\varphi$ then $\psi$, <u>while</u> for the given $s^{l_x}$ or, equivalently, for all $s^{l_n}$ in question, $M, s^{l_n} \vDash \varphi \wedge \psi$; that is, $\varphi$ and $\psi$ are given $\vDash \mathrm{Q}^{\mapsto}{}_x, \mathrm{Q}^{\mapsto}{}_y$ such that $\varphi, \psi \notin p^{sx}{}_x$.
- $M, s^{sx}, s^{l_x} \vDash (\mathfrak{D}^{\mathfrak{K}}\varphi \to \mathfrak{D}^{\mathfrak{K}}\psi)$ *iff* $\varphi$ and $\psi$ are given belief states of the given simultaneous moment, viz., $(b^{sx}{}_x, b^{sx}{}_y) \in s^{sx}$ and if whenever $\varphi$ holds true that $\vDash \mathrm{Q}^{\mapsto}(b^{sx}{}_x)$, then $\psi$ holds true as $\vDash \mathrm{Q}^{\mapsto}(b^{sx}{}_y)$, <u>while</u> for the given $s^{l_x}$ or, equivalently, for all $s^{l_n}$ in question, $M, s^{l_n} \vDash \varphi \wedge \psi$; that is, $\varphi$ and $\psi$ are given $\vDash \mathrm{Q}^{\mapsto}{}_x, \mathrm{Q}^{\mapsto}{}_y$ such that $\varphi, \psi \notin p^{sx}{}_x$.

Assume that $\nvDash (\mathfrak{D}^{\mathfrak{K}}\varphi \to \mathfrak{D}^{\mathfrak{K}}\psi)$. We have that either (i) it is *not* true that both $\varphi$ and $\psi$ are belief states, $(b^{sx}{}_x, b^{sx}{}_y) \in s^{sx}$, (ii) it is *not* true that whenever $\varphi$ holds true that $\vDash \mathrm{Q}^{\mapsto}(b^{sx}{}_x)$, then $\psi$ holds true as $\vDash \mathrm{Q}^{\mapsto}(b^{sx}{}_y)$, (iii) it is *not* true that $\varphi$ and $\psi$ are given $\vDash \mathrm{Q}^{\mapsto}{}_x, \mathrm{Q}^{\mapsto}{}_y$ such that $\varphi, \psi \notin p^{sx}{}_x$ or some combination thereof. It is possible then that $\vDash \mathfrak{D}^{\mathfrak{K}}(\varphi \to \psi)$, as it is possible that (i) $\varphi$ and $\psi$ are pre-belief state moments and whenever $\varphi$ then $\psi$ and (ii) $\varphi$ and $\psi$ are given $\vDash \mathsf{Q}^{\mapsto}{}_{x}, \mathsf{Q}^{\mapsto}{}_{y}$ such that $\varphi, \psi \notin p^{sx}{}_{x}$.

Thus, we do *not* have,

$$\mathfrak{D}^{\mathfrak{K}}(\varphi \rightarrow \psi) \rightarrow (\mathfrak{D}^{\mathfrak{K}}\varphi \rightarrow \mathfrak{D}^{\mathfrak{K}}\psi)$$

■

Then, for the frame condition for **T**,

$$K\varphi \rightarrow \varphi$$

$$\mathfrak{D}^{\mathfrak{K}}\varphi \rightarrow \varphi$$

- Given that if $\vDash \mathfrak{D}^{\mathfrak{K}}\varphi$ then $\vDash \mathsf{Q}^{\mapsto}(b^{sx}{}_{x})$, it is clear that $\mathfrak{D}^{\mathfrak{K}}\varphi \longrightarrow \varphi$. Thus, we have,

$$\mathfrak{D}^{\mathfrak{K}}\varphi \rightarrow \varphi$$

■

### 4.2 Non-normality

We can see that PQG logic—when its semantics are contrasted with those of normal modal logics—is a non-normal logic. We will now examine some non-normal axioms, to get an even better idea of the extreme uniqueness of the logic. In [7] we see two interesting, non-normal axioms, **M** and **C**. We analyze the former of these in terms of PQG below.

For the frame condition for **M**,

$$\Box(\varphi \wedge \psi) \rightarrow \Box\varphi \wedge \Box\psi$$

$$\mathfrak{D}^{\mathfrak{K}}(\varphi \wedge \psi) \rightarrow \mathfrak{D}^{\mathfrak{K}}\varphi \wedge \mathfrak{D}^{\mathfrak{K}}\psi$$

- $M, s^{Sx}, s^{l_x} \vDash \mathfrak{D}^{\mathfrak{K}}(\varphi \wedge \psi)$ *iff* $\varphi$ and $\psi$ are true at given pre-belief state moments.
- $M, s^{Sx}, s^{l_x} \vDash \mathfrak{D}^{\mathfrak{K}}\varphi \wedge \mathfrak{D}^{\mathfrak{K}}\psi$ *iff* $\varphi$ and $\psi$ are given belief states of the given simultaneous moment, viz., $(b^{sx}{}_{x}, b^{sx}{}_{y}) \in s^{Sx}$ and $\varphi$ holds true that $\vDash \mathsf{Q}^{\mapsto}(b^{sx}{}_{x})$, and $\psi$ holds true as $\vDash \mathsf{Q}^{\mapsto}(b^{sx}{}_{y})$.

Assume that $\nvDash \mathfrak{D}^{\mathfrak{K}}\varphi \wedge \mathfrak{D}^{\mathfrak{K}}\psi$. We have either (i) it is *not* true that *both* $\varphi$ and $\psi$ are given belief states of the given simultaneous moment, $(b^{sx}{}_{x}, b^{sx}{}_{y}) \in s^{Sx}$, or (ii) it is *not* true that *both* $\varphi$ holds true that $\vDash \mathsf{Q}^{\mapsto}(\, b^{sx}{}_{x})$, and $\psi$ holds true as $\vDash \mathsf{Q}^{\mapsto}(\, b^{sx}{}_{y})$. According to (i) and (ii), there is no reason why we cannot have that $\vDash \mathfrak{D}^{\mathfrak{K}}(\varphi \wedge \psi)$. Thus, we do *not* have,

$$\mathfrak{D}^{\mathfrak{K}}(\varphi \wedge \psi) \rightarrow \mathfrak{D}^{\mathfrak{K}}\varphi \wedge \mathfrak{D}^{\mathfrak{K}}\psi$$

■

By now, it should be apparent that PQG logic is no usual "modal" logic. Although PQG logic incorporates modalities, it represents a sharp diversion from the usual modal logics, taking what can best be called a *descriptive turn*. Above we have presented the nature of the divergence from the norm. Below, we will now present some principles of PQG.

### 4.3 Discussion of Axioms and Principles

Bringing in the usual alethic modal operators, ◇ and □ (in tandem with our ◎ and ▣), neglecting to expand upon their corresponding frame condition, we can now exhibit the expressivity of PQG. We will start with a brief, non-exhaustive list of principles that are validated by the given semantics, followed by an in-depth discussion of the nature of potential axioms of PQG.

#### 4.3.1 Principles

$$▣p \vDash (\mathfrak{D}^{\mathfrak{B}}p \vee \mathfrak{D}^{\mathfrak{K}}p)$$

$$◎p \vDash \neg(\mathfrak{D}^{\mathfrak{B}}p \vee \mathfrak{D}^{\mathfrak{K}}p)$$

$$(\mathfrak{D}^{\mathfrak{B}}p \vee \mathfrak{D}^{\mathfrak{K}}p) \vDash ▣p$$

$$▣p \vDash [\mathfrak{B}p \wedge \forall\, d.\, d \in \mathbf{\mathfrak{D}}_{fv}{}^{s_{lx}}\, [d \in \mathbb{R}_{fv}]$$

$$\mathfrak{D}^{\mathfrak{B}\mathfrak{M}^{n}}p \vDash \mathfrak{D}^{\mathfrak{B}\mathfrak{M}^{n\text{-}1}}p, \mathfrak{D}^{\mathfrak{B}\mathfrak{M}^{n\text{-}2}}p, \ldots, \mathfrak{D}^{\mathfrak{B}}p$$

$$\mathfrak{D}^{\mathfrak{K}\mathfrak{M}^{n}}p \vDash \mathfrak{D}^{\mathfrak{K}\mathfrak{M}^{n\text{-}1}}p, \mathfrak{D}^{\mathfrak{K}\mathfrak{M}^{n\text{-}2}}p, \ldots, \mathfrak{D}^{\mathfrak{K}}p$$

$$\Box\mathfrak{D}^{\mathfrak{B}} \vDash \Box\mathbf{N}[(b^{sx} ✧ f_{v}{}^{s_{x}}) \wedge ([(s^{l}, s^{s})_{1}, \ldots, (s^{l}, s^{s})_{n}] \between [f_{v}{}^{sx}{}_{1}, \ldots, f_{v}{}^{sx}{}_{n}])]$$

### 4.3.2 Moving towards Axioms

Let us reintroduce a distinction we made in the very beginning of this section. We assert a distinction between (i) the functional assembly of cognition; the entire logically analyzable apparatus from perception and memory to reasoning and qualia and (ii) an apparatus of implication, with all of the accompanying functional assembly abstracted away. Where we are dealing with a logic in the latter sense, its axioms simply "structural" starting points, notions upon which we can analyze reasoning. While, where we are dealing with a logic in the former sense, axioms quite literal build the *cognitive architecture* of the agent in question. Granted that our epistemic/doxastic logic is not dependent upon the traditional modal frame conditions (it eschews these entirely), we cannot rely upon the usual means of axiom consideration and generation. The process of axiom generation for PQG logic will be an in-depth specification of a given agent(s) cognitive architecture, as it is constrained into the PQG template of *volitionary functions*, *belief states* and *psychological concepts*.

As bizarre as this approach may seem, it is not completely insular. There are similarities between PQG and Goertzel's approach in his [11]. Although he does not build upon the quanta of percept, qualia and cognition as we do, he does utilize the same method, building upon the atomic units of *cognits*, *actions*, *goals*, *observations* and *rewards*, which he denotes as,

$$c_1 a_1 o_1 g_1 r_1 c_2 a_2 o_2 g_2 r_2 \ldots$$

Goertzel takes the agent to be a function $\pi$ taking the "current history" as input, producing an action as output. Although Goertzel utilizes hypergraphs and categorical semantics, we can see many similarities between his basic framework and our PQG logic, especially in keeping in mind the former half of the distinction we just mentioned. Until now, it has simply been the case that such psychological-functionalist approaches have been demonstrated in regions other than logic, namely in artificial intelligence. Though the initial comparison of PQG with traditional modal logics may leave one underwhelmed in terms of expressivity, the potential applications into various domains, i.e., the translatability of the system is profound. For instance, in the domain of multi-agent systems, the BDI (belief-desire-intention) software has dominated the field essentially from its inception, and calls for transcending the approach have only just recently arisen [22]. The BDI approach is a technical offshoot from the *philosophical* groundwork laid by Bratman in his [5]; it is crucial here to notice the importance of philosophical, conceptual structures upon which utile technical innovations can arise. Revisiting Schwitzgebel's [30], we can utilize his distinction between *deep* vs. *superficial* accounts of psychological states. We examine two excerpts. First, we examine his main demarcation between the deep and superficial accounts:

> Let's say that relative to a *class of surface phenomena*, an account of a property is *deep* if it identifies possession of the property with some feature other than patterns in those same surface phenomena – some feature that presumably explains or causes or underwrites those surface patterns. In contrast, let's say that an account is superficial if it identifies possession of the property simply with patterns in the surface phenomena. (emphasis added)

In terms of modal logic, we can adumbrate a class of surface phenomena for belief and knowledge, namely the class of all sentences $K\varphi$ and $B\varphi$ and the specified accessibility relation(s). It is quite clear that this is *surface phenomena* of the psychological states of knowledge and belief. In this case, the said "patterns in the surface phenomena" are the various accompanying logical axioms and theorems. In the case of PQG, we aspire to give a "deeper" account of these surface phenomena; we wish to investigate *why* certain sentences are doxastically and epistemically accessible. Enter Schwitzgebel once more.

> Accounts of psychological properties can likewise be deep or superficial *relative* to a class of surface phenomena […] Any account of a psychological property that identifies possession of that property with having a particular folk-psychologically *non-obvious functional architecture* will be deep relative to any class of surface phenomena that does not include folk-psychologically non-obvious functional architecture. In both of these respects, my approach to the attitudes is superficial rather than deep. (emphasis added)

Relative to the class of surface phenomena of all sentences $K\varphi$ and $B\varphi$ and the specified accessibility relation(s) that we named, we can have deep or superficial accounts. In terms of properties of the human psychological "functional architecture", it can certainly be said to be *folk-psychologically obvious* that the process of reasoning is sequential or "programmatic" in nature. This much renders such approaches as seen in [4, 13], in which dynamic logic—the logic of computational, "dynamic" *procedure*—is utilized to give an ostensibly "deeper" account of the class of doxastic and epistemic surface phenomena, relatively shallow for our purposes. Where the former approach provides, in accordance with Schwitzgebel's analysis, a functionally superficial account of the phenomena, we assert that PQG logic—with its novel doxastic and epistemic semantics—serves as a *deep* account of said class of phenomena. With this being said, we may now commence an investigation into potential axioms and classes thereof.

We shall first begin by iterating the bulk of the formal assembly of PQG.

**Distinguished Functions, Sets and Variables**

| | |
|---|---|
| Belief state | ${b^{sx}}_x \in s^s$; ${b^{sx}}_x$ as $Q^{\mapsto}$ at $s^l$ |
| Taking function | $\mathcal{T}$: ${Q^{\mapsto}}_x \mapsto {Q^{\mapsto}}_y$, st. ${Q^{\mapsto}}_x >_l {Q^{\mapsto}}_y$ |
| Forming function | $\mathcal{F}$: ${Q^{\mapsto}}_y \mapsto {Q^{\mapsto}}_z$, st. ${Q^{\mapsto}}_y$ is from a given $\mathcal{T}$ |
| Psychological concept | a *given* mapping of the function $\mathcal{F}$, ${Q^{\mapsto}}_y \mapsto {Q^{\mapsto}}_z$, $c_x$ |
| Volitionary functions | an $n$-tuple, $f_v$, $f_i$, …, $f_n$ of ${\mathfrak{F}_s}^s$ |
| Recommended quanta strings | where each ${Q^{\mapsto}}_z$ is from a forming function, $\mathcal{F}$: ${Q^{\mapsto}}_y \mapsto {Q^{\mapsto}}_z$ and where the output of ${f_x}^n$, $Q^{\mapsto}$ is RQS |
| Pre-belief state moments | ${p^{sx}}_x \in \mathbb{P}$, such that members of $\mathbb{P}$ are members of $\mathbb{B}$ st. $b^{sx}$ ✧ ${f_v}^{sx}$ and $[(s^l, s^s)_1, \ldots, (s^l, s^s)_n]$ ⟠ $[{{f_v}^{sx}}_1, \ldots, {{f_v}^{sx}}_n]$ |
| Prime volitionary rules<br><br><br><br><br>Prime volitionary rules (contd.) | a set $\mathbb{R}_{f_v}$ for the given prime volitionary function $f_v$,<br>$f_i^{\,1}$: $(P^i{}_2\langle(c_i), (Q^{\mapsto}{}_z{}^i)\rangle, P^i{}_4\langle(c_j), (Q^{\mapsto}{}_z{}^j)\rangle) \mapsto P^j{}_4(Q^{\mapsto})$<br>$f_j^{\,2}$: $(P^i{}_6\langle(c_i), (Q^{\mapsto}{}_z{}^i)\rangle, P^i{}_2\langle(c_j), (Q^{\mapsto}{}_z{}^j)\rangle) \mapsto P^j{}_3(Q^{\mapsto})$,<br>…<br>$f_n^{\,n\text{-}1}$: … |
| Volitionary determination | $\mathbb{D}_{f_v}{}^{b_x}$; these are seen to be rules dictating the “ability” of $f_v{}^{s_x}$ to “reach” $b^{sx}{}_x$, st., if these rules hold, per all given $p^{sx}{}_x$, then $b^{sx}{}_x$ will eventually hold; also, a subset, a set $\mathfrak{S}_{f_v}{}^{b_x}$ of *minimal rules* and a *superset*, a set $\mathfrak{X}_{f_v}{}^{b_x}$ of *maximal rules* |

**Distinguished Relations**

| | |
|---|---|
| Volitionary acceptance | $\forall f\,(f \in \mathbb{D}_{f_v}{}^{b_2})\,[f \in \mathbb{R}_{f_v}] \equiv b^{sx} \diamond f_v{}^{s_x}$ |
| `Volitionary invariance` | $\forall\, s\, [(s^l, s^s)_1, \ldots, (s^l, s^s)_n]\, (\forall\, f\, (f \in \mathbb{D}_{f_v}{}^{b_2})\, [f \in \mathbb{R}_{f_v}]) \equiv [(s^l, s^s)_1, \ldots, (s^l, s^s)_n] \mathbin{\text{⟠}} [f_v{}^{s_x}{}_1, \ldots, f_v{}^{s_x}{}_n]$ |
| `Linear ordering relation` | $s_x <_l s_y$, only if $s_x$, $s_y$ are both of sort `lin` |
| `Metasimultaneous ordering relation` | $s_x <_S s_y$, only if $s_x$, $s_y$ are both of sort `sim` |
| `Metasimultaneous inclusion relation` | $s_x \supset_S s_y$, only if $s_x$ is of sort `lin` and $s_y$ is of sort `sim` |
| `Hypothetical ordering relation` | $s_x <_h{}^{b s_x} s_y$, only if $s_x$ and $s_y$ are of sort `pre` |

### 4.3.3 The Structure of Potential Axioms

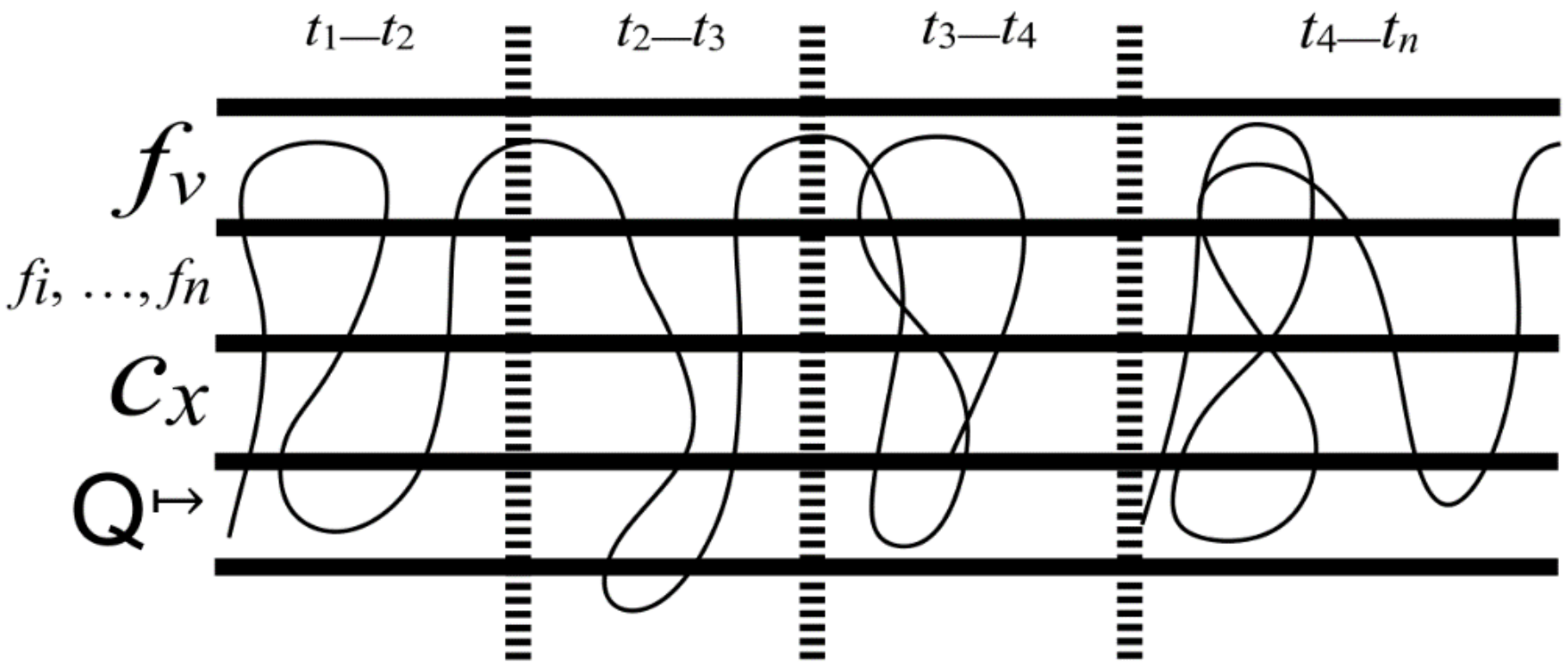

**Fig. 2** Here we see a delineation of the key distinguished classes of functions, $f_v$, the prime volitionary function and a set $\{f_i, \ldots, f_n\}$ of functions whose outputs serve as arguments of $f_v$; we also see two distinguished sets which the aforementioned functions take as their arguments. The zig-zagging of the line represents the confluence of potential domains of axiom consideration available per each “window” of time, $t_x$–$t_y$.

We can start by reintroducing our hypothetical set $\mathbb{D}_{f_v}{}^{b_2}$, through which we gave an example of the various potential predicates accompanying the aspects of the ordered sets, $\langle c_x, \mathrm{Q}^{\mapsto}{}_z{}^x \rangle$, which comprise the arguments of functions $f_i, \ldots, f_n$.

---

Below we list a hypothetical set $\mathbb{D}_{f_v}{}^{b_2}$ consisting of only two functions "below" the prime volitionary function $f_v$,

$$f_i{}^1: (P^i{}_2\langle(c_i), (\mathrm{Q}^{\mapsto}{}_z{}^i)\rangle, P^i{}_4\langle(c_j), (\mathrm{Q}^{\mapsto}{}_z{}^j)\rangle) \;\mapsto\; P^j{}_4(\mathrm{Q}^{\mapsto})$$

$$f_j{}^2: (P^i{}_6\langle(c_i), (\mathrm{Q}^{\mapsto}{}_z{}^i)\rangle, P^i{}_2\langle(c_j), (\mathrm{Q}^{\mapsto}{}_z{}^j)\rangle) \;\mapsto\; P^j{}_3(\mathrm{Q}^{\mapsto}),$$

where the $P^x{}_x$s are given predicates, whose definitions may be of various sorts. We provide example definitions below,

$$P^i{}_6\langle(c_i), (\mathrm{Q}^{\mapsto}{}_z{}^i)\rangle \equiv c_i \ni \forall\, \mathrm{Q}^{\mapsto}{}_y\, (\mathrm{Q}^{\mapsto}{}_y \in \mathcal{F})\, [\mathrm{Q}^{\mapsto}{}_y \in \mathrm{RQS}_{s_x}{}^{f_i}] \ldots \mathrm{Q}^{\mapsto}{}_z{}^i \ni \exists\, g, q\, ((g, q) \in \mathcal{T}\!:\; \mathrm{Q}^{\mapsto}{}_y >_l \mathrm{Q}^{\mapsto}{}_x{}^2$$

$$P^j{}_3(\mathrm{Q}^{\mapsto}) \equiv\; \mathrm{Q}^{\mapsto} \ni \exists\, s^l\, (s^l <_l \mathrm{Q}^{\mapsto})$$

---

Per a given volitionary function of a given agent, we may wish to specify various varieties of predicate, per predicate "class", i.e., kinds of $P^i$ and kinds of $P^j$. From here, as was presented in the most recent figure above, we can specify these predicates—at the broadest level—in terms of the nature of previous $f_v$, $\{f_i, \ldots, f_n\}$, $c_x$ and $\mathrm{Q}^{\mapsto}$. Following the intuitive division of volitionary functions into distinctive "units", as we stated above, citing [xx] as guiding motivation, we will appropriately delineate the set $\{f_i, \ldots, f_n\}$ into various varieties. Per function, we may specify the nature of its arguments, i.e., concept-quanta pairings, in tandem with the nature of its outputs, i.e., recommended quanta strings. Then, per psychological concept ($c_x$), we may specify the nature of its accompanying $\mathcal{T}$ and $\mathcal{F}$, that is, the functions $\mathcal{T}: \mathrm{Q}^{\mapsto}{}_x \mapsto \mathrm{Q}^{\mapsto}{}_y$ and $\mathcal{F}: \mathrm{Q}^{\mapsto}{}_y \mapsto \mathrm{Q}^{\mapsto}{}_z$. And, finally, per quanta string, we can of course specify its nature in various ways.

All of this above leads us to the conclusion that a given agent $\alpha$ has a collection of axioms, $\boldsymbol{T}$, that is delineated into four basic, interdefinable "sets":

| |
|---|
| Set $F_v$ of $f_v$ at $s^s{}_n$ |
| Set $F_n$ of $f_n$ at $s^s{}_n$ |
| Set $C$ of $c_x$ at $s^l{}_n$ |
| Set $Q$ of $\mathsf{Q}^{\mapsto}$ at $s^l{}_n$ |

Although we will forgo a deep dissection of these potential axiom sets, it should be apparent how various possible schemas can be constructed.

## 5 On Logical Omniscience

We mentioned in the opening section that there have been many recent, innovative approaches to solving the problem of logical omniscience [4, 13, 28]. We also mentioned that *all* of these approaches utilize Kripke semantics (whether through dynamic logic [4, 13], “possibilist” quantification and novel operators [28], etc.). However, there is something additional that we should mention. For all the myriad standpoints in the philosophy of belief, in contemporary logical systems, only one approach has been taken—what I shall term *sententialism*. The grammatical nature of belief, as explored marvelously in [15], seems to have taken deep root in the logics of belief. In other words, it proves difficult, in many doxastic logics, to find where “George believes that capitalism is terrific” ends and where $B_G(t(c))$ begins. Though there is assuredly an abundance of reasons why this seems to be the case, one core reason seems to be the lasting impact of Hintikka’s endeavor, as explicated in [14], to pragmatically pursue the question of what it *means* to know, while neglecting the question of what it *is* to know. Where we are solely concerned with the former matter, it makes no difference what sort of relation $B_G(t(c))$ has to “George believes that capitalism is terrific”, as we can examine the meaning (the logical implications) of the sentence simply through *direct translation*, from informal sentence to formal expression.

Even after the relevant changes have been implemented and the agent does not succumb to logical omniscience, as in [4, 13, 28], logical expressions such as $B_G(t(c))$ still derive their structure from accompanying informal sentences. In PQG logic, we see a fundamental shift from this sententialism. We no longer care about the complement clause, “that capitalism is terrific”; we care about *descriptive* matters. Thus, an example of an informal mold for a PQG expression would instead be something like, “George’s belief is a bundle of qualia and cognitions”, which, of course, is not very informative on its face. However, once it is integrated into the cognitive-functional apparatus of PQG (volitionary functions, concepts, etc.), it becomes extremely useful. Keeping with its descriptive tone, our semantics for such a *sententialist* sentence as “George believes that if it is raining then the government has made it so” adjust it to a descriptively accommodating “If some hypothetical precursor to George’s belief that it is raining is met, then some hypothetical precursor to his belief that the government is behind something is met”. It is in this manner that we avoid integrating the content of the complement clause into the logical implications of the belief itself. This is because—instead of dealing with the *meaning* of the *logical implications* of that which is believed—we are dealing with the nature of the very *beliefs themselves* and their necessary *functional criteria*.

What implications does this have for the problem of logical omniscience then? Most broadly, the problem can be seen as a pesky feature of possible worlds semantics; it is built into it the very framework. We may defer to Bjerring and Skipper, from their [4]:

> To see why [logical omniscience is an issue], suppose that an agent believes a proposition *p*, and let *q* be any logical consequence of *p*. Since the agent believes *p*, *p* is true at all possible worlds that are doxastically possible for the agent. And since *p* entails *q*, all possible worlds that verify *p* also verify *q*. Hence *q* is true at all doxastically possible worlds for the agent, which means that the agent believes *q*. So if the agent believes *p*, she believes all logical consequences of *p*.

More specifically then, for all normal modal logics, there are accompanying principles of *closure*: (i) closure under known implication, (ii) closure under conjunction elimination and (iii) closure under disjunction introduction. As we have already shown in our comparison of PQG with the frame condition of **K**, our logic is non-normal. Moreover, given that its semantics of knowledge and belief are outside of possible-worlds semantics altogether, it is not “non-normal” in the traditional Kripkean sense of the term. Given that sentences of the form $\mathfrak{D}^{\mathfrak{K}}(\ldots\varphi\ldots)$ are semantically defined on the basis of *pre-belief state* moments, that is—hypothetical *bundles* of quanta leading up to a belief as quanta—we can effectively side-step the issues associated with modal distribution, i.e., the binding of a generic modal operator $\bigcirc$ with a non-atomic WFF, say, $\varphi \wedge \psi$: $\bigcirc(\varphi \wedge \psi)$. In deontic logic, such difficulties have long been widespread, as evinced by Chisholm’s paradox. While, in doxastic and epistemic logic, such issues seem to have been emphasized to lesser degree. “X knows that if *x* then *y*” is astoundingly different than “If X knows *x* then he knows that *y*”; this much has, to some degree, of course, been recognized and accounted for. However, the extent to which the endogenous/exogenous reasoning schism—i.e., $B\varphi \rightarrow B\psi$ vs. $B(\varphi \rightarrow \psi)$—has hampered the development of doxastic and epistemic logics, I contend, has been largely neglected. In any case, PQG is well-equipped to tackle all principles of doxastic and epistemic closure (which, by and large, arise from said endogenous/exogenous schism).

For (i), closure under known implication, $(K\varphi \wedge K((\varphi \rightarrow \psi)) \rightarrow K\psi$, we do not have a PQG isomorphic counterpart. However, introducing an operator $\mathfrak{D}^{\mathfrak{P}}$, meaning that the following sentence(s) is a pre-belief state moment, we *do* have, $(\mathfrak{D}^{\mathfrak{K}}\varphi \wedge \mathfrak{D}^{\mathfrak{K}}((\varphi \rightarrow \psi)) \rightarrow \mathfrak{D}^{\mathfrak{P}}\psi$. Granted, this is rather trivial, but the point is that the closure principle is successfully obviated. For the other two closure principles, such sentences can be formulated as well. For (ii) closure under conjunction elimination, we have $\mathfrak{D}^{\mathfrak{K}}(\varphi \wedge \psi) \rightarrow \mathfrak{D}^{\mathfrak{P}}\varphi$ but *not* $\mathfrak{D}^{\mathfrak{K}}(\varphi \wedge \psi) \rightarrow \mathfrak{D}^{\mathfrak{K}}\varphi$. And, for closure under disjunction introduction, we have $\mathfrak{D}^{\mathfrak{P}}\varphi \rightarrow \mathfrak{D}^{\mathfrak{K}}(\varphi \vee \psi)$ but *not* $\mathfrak{D}^{\mathfrak{K}}\varphi \rightarrow \mathfrak{D}^{\mathfrak{K}}(\varphi \vee \psi)$.

We have made it clear that our epistemic and doxastic operators, $\mathfrak{D}^{\mathfrak{K}}$ and $\mathfrak{D}^{\mathfrak{B}}$, are *descriptive* in nature, thus making the comparison between PQG and traditional doxastic and epistemic logics tricky. But, for all the ostensible “flatness” of entailment, we can likely move from $\mathfrak{D}^{\mathfrak{K}}(\ldots\varphi\ldots)$, the assertion of some arrangement of pre-belief state moments, to $\mathfrak{D}^{\mathfrak{K}}\varphi\ldots$, the assertion of some $\varphi\ldots$ being a $\vDash \mathsf{Q}^{\mapsto}(\, b^{sx}{}_{x})$. There is good reason to believe the following principle:

$$(\mathfrak{D}^{\mathfrak{K}}\varphi \wedge \mathfrak{D}^{\mathfrak{K}}((\varphi \equiv \psi)) \rightarrow \mathfrak{D}^{\mathfrak{K}}\psi,$$

which states that if there is some $\psi$ that is “everywhere coupled” with $\varphi$, in terms of the ordering $<_h^{bs_x}$, *and* $\varphi$ is known, then that $\psi$ is known as well. For example, say that we have a $\vDash \mathsf{Q}^{\mapsto}(\, b^{sx}{}_{x})$, interpreted as a quanta sequence corresponding to the fact “it is raining”. Then, given a $\vDash \mathsf{Q}^{\mapsto}(\, b^{sx}{}_{y})$, for the *pre-belief state* moments of this, we have a $\vDash \mathsf{Q}^{\mapsto}(\, p^{sx}{}_{y})$ interpreted as a quanta sequence corresponding to the fact “agent $\alpha$ looks out the window”. Suppose it is the case that whenever we have the quanta sequence of $\mathsf{Q}^{\mapsto}(\, b^{sx}{}_{x})$ as a pre-belief state, $\mathsf{Q}^{\mapsto}(\, p^{sx}{}_{x})$, we also have $\mathsf{Q}^{\mapsto}(\, p^{sx}{}_{y})$; suppose that the converse holds as well. If we simply had one half of the biconditional, it could simply be the case that the “agent $\alpha$ looks out the window” happens to be a low-level concomitant of “it is raining”, i.e., whenever “it is raining”, “agent $\alpha$ looks out the window”, *in the context of a given pre-belief state sequence*, does not entail that the latter corresponding quanta sequence is one of the fulfillment of belief, as the governing belief of the sequence could render the two states to be inconsequentially related. However, when the biconditional holds, there is necessarily some sort of interdependency between “it is raining” and “agent $\alpha$ looks out the window”. And, given that (i) the former is already a belief of agent $\alpha$, (ii) the former shows up as a pre-belief state of another belief and (iii) the former and the latter are pre-belief state-wise interdependent, we can say that the two facts form some sort of belief complex; thus $(\mathfrak{D}^{\mathfrak{K}}\varphi \wedge \mathfrak{D}^{\mathfrak{K}}((\varphi \equiv \psi)) \rightarrow \mathfrak{D}^{\mathfrak{K}}\psi$.

## 6 Summary

We began by noting the lack of contemporary interest in a logical analysis of what belief is, as opposed to what it means, via its implications. Noting the inadequacy of possible-worlds semantics to provide a psychologico-functional account of belief and knowledge, we developed a novel, functionally-sophisticated semantics, terming the accompanying logic “PQG logic”. Taking quanta of percept, qualia and cognition as logical bedrock, we build up to modality via the interplay of quanta-carrying functions, termed volitionary functions. Having reached our satisfaction definitions for our modalities, we discussed potential axioms and the contribution of PQG logic to the problem of logical omniscience.